\documentclass{amsart}
\usepackage{lastpage}
\usepackage{yhmath,verbatim,mathabx}
\usepackage[all]{xy}

\DeclareFontFamily{U}{mathx}{\hyphenchar\font45}
\DeclareFontShape{U}{mathx}{m}{n}{
	<5> <6> <7> <8> <9> <10>
	<10.95> <12> <14.4> <17.28> <20.74> <24.88>
	mathx10
}{}
\DeclareSymbolFont{mathx}{U}{mathx}{m}{n}
\DeclareFontSubstitution{U}{mathx}{m}{n}
\DeclareMathAccent{\widecheck}{0}{mathx}{"71}
\DeclareMathAccent{\wideparen}{0}{mathx}{"75}

\newcommand{\bdis}{\begin{displaymath}}
\newcommand{\edis}{\end{displaymath}}
\newcommand{\be}{\begin{equation}}
\newcommand{\ee}{\end{equation}}
\newcommand{\mbb}{\mathbb}
\newcommand{\mcal}{\mathcal}

\newcommand{\vp}{\varphi}

\newcommand{\zf}{\zeta\left(\frac{1}{2}+it\right)}

\DeclareMathOperator{\re}{Re} 
\DeclareMathOperator{\im}{Im}

\theoremstyle{definition}

\theoremstyle{remark}
\newtheorem{remark}[]{Remark}

\newtheorem*{mydef11}{{\bf Theorem 1}}

\newtheorem*{mydef12}{{\bf Theorem 2}}

\newtheorem*{mydef41}{{\bf Corollary 1}}

\newtheorem*{mydef42}{{\bf Corollary 2}}

\newtheorem*{mydef43}{{\bf Corollary 3}}

\newtheorem*{mydef44}{{\bf Corollary 4}}

\newtheorem*{mydef51}{{\bf Lemma 1}}

\newtheorem*{mydef52}{{\bf Lemma 2}}

\newtheorem*{mydef53}{{\bf Lemma 3}}

\newtheorem*{mydef54}{{\bf Lemma 4}}

\newtheorem*{mydefECHF1}{{\bf Exact Complete Hybrid Formula 1}} 

\newtheorem*{mydefECHF2}{{\bf Exact Complete Hybrid Formula 2}} 

\newtheorem*{mydefACHF1}{{\bf Asymptotic Complete Hybrid Formula 1}} 

\newtheorem*{mydefACHF2}{{\bf Asymptotic Complete Hybrid Formula 2}}

\numberwithin{equation}{section}

\begin{document}

\title{Jacob's ladders and grafting of the complete hybrid formulas into $\zeta$-synergetic meta-functional equation for Riemann's zeta-function} 

\author{Jan Moser}

\address{Department of Mathematical Analysis and Numerical Mathematics, Comenius University, Mlynska Dolina M105, 842 48 Bratislava, SLOVAKIA}

\email{jan.mozer@fmph.uniba.sk}


\begin{abstract}
In this paper we prove that there is a grafting of complete hybrid formula on meta-functional equation. This synergetic formula gives the result of a $\zeta$-chemical reaction between different sets of values of Riemann's zeta-function on the critical strip. 

\centering Dedicated to recalling on A.S. Eddington's - manipulations with set of world-constants. 
\end{abstract}
\maketitle

\section{Introduction}

In this paper we obtain the first set of meta-functional equations. It is an example of entirely new kind of formulas in the theory of Riemann's zeta-function. Namely: there are infinite sets 
\be \label{1.1} 
\begin{split}
& \{\alpha_1^{1,1}\},\ \{\alpha_1^{2,1}\},\ \{\alpha_1^{3,1}\}\in (L,+\infty),\ L>0, \\ 
& \{ w_1\},\ \{ w_2\},\ \{ w_3\}\in S_{\sigma_1,\sigma_2}=(\sigma_1,\sigma_2)\times (0,+\infty), \\ 
& \frac 12<\sigma_1<\sigma_2<1, \\ 
& w_1,w_2,w_3\in\mbb{C} 
\end{split}
\ee 
of values such that the following meta-functional equation 
\be \label{1.2} 
\begin{split}
& \left|\zeta\left(\frac 12+i\alpha_1^{2,1}\right)\right|^2|\zeta(w_2)|-
\left|\zeta\left(\frac 12+i\alpha_1^{3,1}\right)\right|^2|\zeta(w_1)|\sim \\ 
& \sim \left|\zeta\left(\frac 12+i\alpha_1^{3,1}\right)\right|^2|\zeta(w_3)|,\ L\to\infty 
\end{split}
\ee 
holds true. 
 
\begin{remark}
The formula (\ref{1.2}) is: 
\begin{itemize}
\item[(a)] the simplest element of meta-functional equations obtained in this paper, 
\item[(b)] synergetic one since it is a result of cooperative interactions (see our interpretation in \cite{8}) between the following sets of values 
\be \label{1.3} 
\left\{\left|\zf\right|^2\right\},\ t\in (L,+\infty),\ |\zeta(w)|,\ w\in S_{\sigma_1,\sigma_2}
\ee 
where, of course, 
\be \label{1.4} 
\left\{\frac 12+it:\ t>L\right\}\bigcap S_{\sigma_1,\sigma_2}=\emptyset. 
\ee 
\end{itemize}   
\end{remark} 

\begin{remark}
Let us notice explicitly that in the case (\ref{1.3}) we have synergetic (cooperative) interactions between the following parts of the function 
\bdis 
\zeta(s),\ s\in\mbb{C}\setminus \{ 1\},\ s=\sigma+it , 
\edis  
namely, between the part 
\bdis 
\left\{\left|\zf\right|^2\right\},\ t>L 
\edis  
on the critical line $\sigma=\frac 12$, and the part (see (\ref{1.1})) 
\bdis 
\{|\zeta(\sigma+it)|\},\ \sigma,\ t\in S_{\sigma_1,\sigma_2} 
\edis  
on the strip (see also (\ref{1.4})). 
\end{remark}

\subsection{} 

Next, we give the following 

\begin{remark}
The proof of meta-functional equation (comp. (\ref{1.2})) is based: 
\begin{itemize}
	\item[(a)] on new notions and methods in the theory of Riemann's zeta-function we have introduced in our series of 46 papers concerning Jacob's ladders (these can be found in arXiv [math.CA] starting with the paper \cite{1}), 
	\item[(b)] on the classical H. Bohr's theorem in analysis (1914) about the roots of the equation 
	\bdis 
	\zeta(s)=a,\ a\not=0,\ a\in\mbb{C},\ \frac 12<\sigma_1<\sigma<\sigma_2<1, 
	\edis  
	\item[(c)] and on our presently introduced notion of \emph{gifting} of complete hybrid formulas that serves for synthesis of the conceptions (a) and (b). 
\end{itemize}
\end{remark}

From our notions and methods (comp. (a)) we use here especially: Jacob's ladder (see \cite{1}), $\zeta$-disconnected set (see \cite{3}), algorithm for generating the $\zeta$-factorization formulas (see \cite{4}) crossbreeding, secondary crossbreeding, exact and asymptotic complete hybrid formula (see \cite{6}, \cite{8}). Short surveys of these notions can be found in \cite{5} and \cite{8}. 

\begin{remark}
We notice explicitly that the notion of meta-functional equation introduced here constitutes also  
the generic complement to the Riemann's functional equation on the critical strip (comp. also \cite{7}). 
\end{remark} 

\section{Zeta-factorization formulas} 

Here we obtain, by making use of our algorithm for generating $\zeta$-factorization formulas (see \cite{5}, (3.1) -- (3.11), comp. \cite{4}) the following. 

\subsection{}  

First, we have the following lemmas (comp. \cite{6}, (2.1) -- (2.9)): 

Since 
\bdis 
\frac 1U\int_{\pi L}^{\pi L+U}\sin^2t{\rm d}t=\frac 12\left(1-\frac{\sin 2U}{2U}\right) 
\edis  
then we have the following 

\begin{mydef51}
For the function 
\be \label{2.1} 
f_1(t)=\sin^2t\in \tilde{C}_0[\pi L,\pi L+U],\ U\in (0,\pi/2),\ L\in\mbb{N} 
\ee  
there are vector-valued functions 
\be \label{2.2} 
\begin{split}
& (\alpha_0^{1,k_1},\alpha_1^{1,k_1},\dots,\alpha_{k_1}^{1,k_1},\beta_1^{k_1},\dots,\beta_{k_1}^{k_1}), \\ 
& 1\leq k_1\leq k_0,\ k_1,k_0\in\mbb{N}
\end{split}
\ee 
(here we fix arbitrary $k_0$) such that the following exact $\zeta$-factorization formula 
\be \label{2.3} 
\prod_{r=1}^{k_1}\frac{\tilde{Z}^2(\alpha_r^{1,k_1})}{\tilde{Z}^2(\beta_r^{k_1})}=
\frac 12\left(1-\frac{\sin 2U}{2U}\right) \frac{1}{\sin^2(\alpha_0^{1,k_1})},\ 
\forall\- L\geq L_0>0 
\ee  
(where $L_0$ is sufficiently big) holds true, where 
\be \label{2.4} 
\begin{split}
& \alpha_r^{1,k_1}=\alpha_r(U,\pi L,k_1;f_1),\ r=0,1,\dots,k_1, \\ 
& \beta_r^{k_1}=\beta_r(U,\pi L,k_1),\ r=1,\dots,k_1, \\ 
& \alpha_0^{1,k_1}\in (\pi L,\pi L+U), \\ 
& \alpha_r^{1,k_1},\beta_r^{k_1}\in (\overset{r}{\wideparen{\pi L}},\overset{r}{\wideparen{\pi L+U}}),\ r=1,\dots,k_1, 
\end{split}
\ee 
and the segment 
\bdis 
(\overset{r}{\wideparen{\pi L}},\overset{r}{\wideparen{\pi L+U}}) 
\edis  
is the $r$-th reverse iteration (by means of the Jacob's ladder, see \cite{3}) of the basic segment 
\bdis 
[\pi L,\pi L+U]=[\overset{0}{\wideparen{\pi L}},\overset{0}{\wideparen{\pi L+U}}]. 
\edis  
Next, 
\be \label{2.5} 
\tilde{Z}^2(t)=\frac{|\zf|^2}{\omega(t)},\ \omega(t)=\left\{1+\mcal{O}\left(\frac{\ln\ln t}{\ln t}\right)\right\}\ln t, 
\ee 
see \cite{2}, (6.1), (6.7), (7.7), (7.8) and (9.1). 
\end{mydef51}

Since 
\bdis 
\frac 1U\int_{\pi L}^{\pi L+U}\cos^2t{\rm d}t=\frac 12\left(1+\frac{\sin 2U}{2U}\right) 
\edis  
then we have the following  

\begin{mydef52}
	For the function 
	\be \label{2.6} 
	f_2(t)=\cos^2t\in \tilde{C}_0[\pi L,\pi L+U],\ U\in (0,\pi/2),\ L\in\mbb{N} 
	\ee  
	there are vector-valued functions 
	\be \label{2.7} 
	\begin{split}
		& (\alpha_0^{2,k_2},\alpha_1^{2,k_2},\dots,\alpha_{k_2}^{2,k_2},\beta_1^{k_2},\dots,\beta_{k_2}^{k_2}), \\ 
		& 1\leq k_2\leq k_0,\ k_2\in\mbb{N}
	\end{split}
	\ee 
	such that the following exact $\zeta$-factorization formula 
	\be \label{2.8} 
	\prod_{r=1}^{k_2}\frac{\tilde{Z}^2(\alpha_r^{2,k_2})}{\tilde{Z}^2(\beta_r^{k_2})}=
	\frac 12\left(1+\frac{\sin 2U}{2U}\right) \frac{1}{\cos^2(\alpha_0^{2,k_2})},\ 
	\forall\- L\geq L_0 
	\ee  
	holds true, where 
	\be \label{2.9} 
	\begin{split}
		& \alpha_r^{2,k_2}=\alpha_r(U,\pi L,k_2;f_2),\ r=0,1,\dots,k_2, \\ 
		& \beta_r^{k_2}=\beta_r(U,\pi L,k_2),\ r=1,\dots,k_2, \\ 
		& \alpha_0^{2,k_2}\in (\pi L,\pi L+U), \\ 
		& \alpha_r^{2,k_2},\beta_r^{k_2}\in (\overset{r}{\wideparen{\pi L}},\overset{r}{\wideparen{\pi L+U}}),\ r=1,\dots,k_2. 
	\end{split}
	\ee 
\end{mydef52} 

\subsection{} 

Secondly, since 
\bdis 
\frac 1U\int_{\pi L}^{\pi L+U}\cos 2t{\rm d}t=\frac{\sin 2U}{2U}
\edis  
then we have the following  

\begin{mydef53}
	For the function 
	\be \label{2.10} 
	f_3(t)=\cos 2t\in \tilde{C}_0[\pi L,\pi L+U],\ U\in (0,\pi/4),\ L\in\mbb{N} 
	\ee  
	there are vector-valued functions 
	\be \label{2.11} 
	\begin{split}
		& (\alpha_0^{3,k_3},\alpha_1^{3,k_3},\dots,\alpha_{k_3}^{3,k_3},\beta_1^{k_3},\dots,\beta_{k_3}^{k_3}), \\ 
		& 1\leq k_3\leq k_0,\ k_3\in\mbb{N}
	\end{split}
	\ee 
	such that the following exact $\zeta$-factorization formula 
	\be \label{2.12} 
	\prod_{r=1}^{k_3}\frac{\tilde{Z}^2(\alpha_r^{3,k_3})}{\tilde{Z}^2(\beta_r^{k_3})}=
	\frac{\sin 2U}{2U} \frac{1}{\cos (2\alpha_0^{3,k_3})},\ 
	\forall\- L\geq L_0 
	\ee  
	holds true, where 
	\be \label{2.13} 
	\begin{split}
		& \alpha_r^{3,k_3}=\alpha_r(U,\pi L,k_3;f_3),\ r=0,1,\dots,k_3, \\ 
		& \beta_r^{k_3}=\beta_r(U,\pi L,k_3),\ r=1,\dots,k_3, \\ 
		& \alpha_0^{2,k_3}\in (\pi L,\pi L+U), \\ 
		& \alpha_r^{3,k_3},\beta_r^{k_3}\in (\overset{r}{\wideparen{\pi L}},\overset{r}{\wideparen{\pi L+U}}),\ r=1,\dots,k_3. 
	\end{split}
	\ee 
\end{mydef53}  

\begin{remark}
Let us remind that 
\bdis 
\begin{split}
& f(t)\in \tilde{C}_0[\pi L,\pi L+U] \ \Leftrightarrow  \\ 
&  f(t)\in C[\pi L,\pi L+U]\wedge f(t)\geq 0,\notequiv 0 \wedge U=o\left(\frac{L}{\ln L}\right),\ L\to\infty. 
\end{split}
\edis 
\end{remark} 

\section{The set of complete hybrid formulas} 

\subsection{} 

Crossbreeding in the set 
\be \label{3.1} 
\{(2.3),(2.8)\},\ U\in (0,\pi/4) 
\ee
of two $\zeta$-factorization formulas (one stage is sufficient) gives the following (comp. \cite{5}) 

\begin{mydefECHF1}
\be \label{3.2}
\begin{split}
& \left\{
\prod_{r=1}^{k_1}\frac{\tilde{Z}^2(\alpha_r^{1,k_1})}{\tilde{Z}^2(\beta_r^{k_1})}
\right\}\sin^2(\alpha_0^{1,k_1})+
\left\{
\prod_{r=1}^{k_2}\frac{\tilde{Z}^2(\alpha_r^{2,k_2})}{\tilde{Z}^2(\beta_r^{k_2})}
\right\}\cos^2(\alpha_0^{2,k_2})=1,\\ 
& \forall\- L\geq L_0,\ 1\leq k_1,k_2\leq k_0,\ U\in(0,\pi/4). 
\end{split}
\ee 
\end{mydefECHF1} 

If we use (\ref{2.7}) and a little of algebra (comp. \cite{8}, subsection 8.2) in (\ref{3.2}) then we obtain the following 

\begin{mydefACHF1}
	\be \label{3.3}
	\begin{split}
		& \left\{
		\prod_{r=1}^{k_1}\frac{\left|\zeta(\frac 12+i\alpha_r^{1,k_1})\right|^2}{\left|\zeta(\frac 12+i\beta_r^{k_1})\right|^2}
		\right\}\sin^2(\alpha_0^{1,k_1})+
		\left\{
		\prod_{r=1}^{k_2}\frac{\left|\zeta(\frac 12+i\alpha_r^{2,k_2})\right|^2}{\left|\zeta(\frac 12+i\beta_r^{k_2})\right|^2}
		\right\}\cos^2(\alpha_0^{2,k_2})\\ 
		& \sim 1,\ L\to\infty . 
	\end{split}
	\ee 
\end{mydefACHF1}  

\subsection{} 

Next, we make the following crossbreeding in the set 
\be \label{3.4} 
\{(2.3),(2.8),(2.12)\},\ U\in(0,\pi/4) 
\ee 
of three $\zeta$-factorization formulas. Namely: 
\begin{itemize}
	\item[(a)] first (\ref{2.8})$\ominus$(\ref{2.3})$\Rightarrow$ 
\be \label{3.5} 
\begin{split}
& \left\{
\prod_{r=1}^{k_2}\frac{\tilde{Z}^2(\alpha_r^{2,k_2})}{\tilde{Z}^2(\beta_r^{k_2})}
\right\}\cos^2(\alpha_0^{2,k_2})-
\left\{
\prod_{r=1}^{k_1}\frac{\tilde{Z}^2(\alpha_r^{1,k_1})}{\tilde{Z}^2(\beta_r^{k_1})}
\right\}\sin^2(\alpha_0^{1,k_1})=\\ 
& =\frac{\sin 2U}{2U}, 
\end{split}
\ee 
\item[(b)] second, the formula (\ref{3.5}) together with the formula (\ref{2.12}) imply the following
\end{itemize} 

\begin{mydefECHF2}
\be \label{3.6}
\begin{split}
	& \left\{
	\prod_{r=1}^{k_2}\frac{\tilde{Z}^2(\alpha_r^{2,k_2})}{\tilde{Z}^2(\beta_r^{k_2})}
	\right\}\cos^2(\alpha_0^{2,k_2})-
	\left\{
	\prod_{r=1}^{k_1}\frac{\tilde{Z}^2(\alpha_r^{1,k_1})}{\tilde{Z}^2(\beta_r^{k_1})}
	\right\}\sin^2(\alpha_0^{1,k_1})=\\ 
	& =\left\{
	\prod_{r=1}^{k_3}\frac{\tilde{Z}^2(\alpha_r^{3,k_3})}{\tilde{Z}^2(\beta_r^{k_3})}
	\right\}\cos(2\alpha_0^{3,k_3}), \\ 
	& \forall\- L\geq L_0,\ 1\leq k_1,k_2,k_3\leq k_0,\ U\in(0,\pi/4). 
\end{split}
\ee 
\end{mydefECHF2} 

If we rewrite (\ref{3.6}) into the form 
\bdis 
\{\dots\}\cos^2(\alpha_0^{2,k_2})=\{\dots\}\sin^2(\alpha_0^{1,k_1})+\{\dots\}\cos(2\alpha_0^{3,k_3})
\edis  
then we obtain, similarly to (\ref{3.3}) the following 

\begin{mydefACHF2}
	\be \label{3.7}
\begin{split}
	& \left\{
	\prod_{r=1}^{k_2}\frac{\left|\zeta(\frac 12+i\alpha_r^{2,k_2})\right|^2}{\left|\zeta(\frac 12+i\beta_r^{k_2})\right|^2}
	\right\}\cos^2(\alpha_0^{2,k_2})-
	\left\{
	\prod_{r=1}^{k_1}\frac{\left|\zeta(\frac 12+i\alpha_r^{1,k_1})\right|^2}{\left|\zeta(\frac 12+i\beta_r^{k_1})\right|^2}
	\right\}\sin^2(\alpha_0^{1,k_1})\\ 
	& \sim \left\{
	\prod_{r=1}^{k_3}\frac{\left|\zeta(\frac 12+i\alpha_r^{3,k_3})\right|^2}{\left|\zeta(\frac 12+i\beta_r^{k_3})\right|^2}
	\right\}\cos(2\alpha_0^{3,k_3}),\ L\to\infty . 
\end{split}
\ee 
\end{mydefACHF2} 

\subsection{} 

In connection with the crossbreeding we give the following 

\begin{remark}
Our description of the operation of crossbreeding (see \cite{6}, comp. \cite{8}, subsection 3.3) 
contains the following expression: 
\begin{center}
	\dots\ that is: after the finite number of eliminations of external functions 
	\bdis 
	E_m(U,T),\ m=1,\dots,M,\dots 
	\edis 
\end{center} 
However, this is not exact. Exact phrase should be as follows: 
\begin{center}
\dots \ that is: after the finite number of eliminations of the variables $U,T$ from the set of external functions \dots 
\end{center} 
We call these $U,T$ as external variables. 
\end{remark}

\section{Theorems about grafting of complete hybrid formulas into meta-functional equations} 

We use the following consequence of the classical H. Bohr's result (see \cite{9}, pp. 76, 99). 

\begin{mydef54}
There is the infinite set of roots of equation 
\be \label{4.1} 
\zeta(s)=a,\ a\not=0,\ a\in\mbb{C}, s=\sigma+it 
\ee 
in the strip 
\be \label{4.2} 
(\sigma_1,\sigma_2)\times (0,+\infty),\ \frac 12<\sigma_1<\sigma_2<1. 
\ee 
\end{mydef54} 

\subsection{} 

The formula (\ref{3.6}) contains the following functions 
\be \label{4.3} 
\sin^2(\alpha_0^{1,k_1}),\ \cos^2(\alpha_0^{2,k_2}),\ \cos(2\alpha_0^{3,k_3})
\ee  
where (see (\ref{2.4}),(\ref{2.9}), (\ref{2.13})) 
\be \label{4.4} 
\begin{split}
& \alpha_0^{1,k_1}=\alpha_0(U,\pi L,k_1;f_1)\in (\pi L,\pi L+U),\ U\in (0,\pi/4), \\ 
& \alpha_0^{2,k_2}=\alpha_0(U,\pi L,k_2;f_2)\in (\pi L,\pi L+U),\ U\in (0,\pi/4),  \\ 
& \alpha_0^{3,k_3}=\alpha_0(U,\pi L,k_3;f_3)\in (\pi L,\pi L+U),\ U\in (0,\pi/4) .  
\end{split}
\ee  
Now, (\ref{4.3}) and (\ref{4.4}) imply 
\be \label{4.5} 
\sin^2(\alpha_0^{1,k_1})\in (0,1/2),\ \cos^2(\alpha_0^{2,k_2})\in (1/2,1),\ \cos(2\alpha_0^{3,k_3})\in(0,1), 
\ee 
i.e. it is true that 
\be \label{4.6} 
\sin^2(\alpha_0^{1,k_1}), \cos^2(\alpha_0^{2,k_2}), \cos(2\alpha_0^{3,k_3}) \in (0,1) 
\ee 
for every 
\bdis 
U\in (0,\pi/4), L\geq L_0, (k_1,k_2,k_3):\ 1\leq k_1,k_2,k_3\leq k_0. 
\edis  

\subsection{} 

Next, we choose (see (\ref{4.1})) 
\bdis 
\sigma_0^1,\sigma_0^2,\sigma_0^3:\ \frac 12<\sigma_1<\sigma_0^1<\sigma_0^2<\sigma_0^3<\sigma_2<1 
\edis  
such that 
\be \label{4.7} 
\begin{split}
& \sigma_1<\sigma_0^1-\delta,\ \sigma_0^3+\delta<\sigma_2, \\ 
& \sigma_0^1+\delta<\sigma_0^2-\delta,\ \sigma_0^2+\delta<\sigma_0^3-\delta 
\end{split}
\ee 
for some sufficiently small $\delta>0$. Now, in the strip 
\be \label{4.8} 
S^1_{\sigma_0^1,\delta}=(\sigma_0^1-\delta,\sigma_0^1+\delta)\times (0,+\infty)
\ee 
there is (see Bohr's Lemma 4) for every fixed and admissible parameters 
\bdis 
U,\pi L,k_1,f_1,\sigma_0^1-\delta,\sigma_0^1+\delta 
\edis 
infinite set 
\be \label{4.9} 
\begin{split}
& W_1=W_1\{\sin^2[\alpha_0(U,\pi L,k_1;f_1)],\sigma_0^1-\delta,\sigma_0^1+\delta\}= \\ 
& = W_1(U,\pi L,k_1,f_1,\sigma_0^1-\delta,\sigma_0^1+\delta)
\end{split} 
\ee 
of elements 
\be \label{4.10} 
w_1=w_1(U,\pi L,k_1,f_1,\sigma_0^1-\delta,\sigma_0^1+\delta)\in W_1\subset 
S^1_{\sigma_0^1,\delta}
\ee 
for which (see (\ref{4.1}),(\ref{4.6})) 
\be \label{4.11} 
\zeta(w_1)=\sin^2(\alpha_0^{1,k_1})>0. 
\ee 
Since (\ref{4.11})$\Rightarrow$ 
\bdis 
\re\{\zeta(w_1)\}=\sin^2(\alpha_0^{1,k_1}),\ \im\{\zeta(w_1)\}=\sin^2(\alpha_0^{1,k_1})=0,
\edis 
then we have the first infinite set of grafts (=equalities) 
\be \label{4.12} 
\sin^2(\alpha_0^{1,k_1})=|\zeta(w_1)|,\ w_1\in W_1
\ee 
generated by the function 
\bdis 
\sin^2(\alpha_0^{1,k_1}) 
\edis  
for every set of fixed and admissible parameters 
\bdis 
U,\pi L,k_1,f_1,\sigma_0^1-\delta,\sigma_0^1+\delta . 
\edis

\subsection{} 

Next, similarly to (\ref{4.7}) -- (\ref{4.12}) we obtain the following objects 
\be \label{4.13} 
S^l_{\sigma_0^l,\delta},W_l\subset S^l_{\sigma_0^l,\delta},\ l=2,3
\ee 
and two other infinite set of grafts 
\be \label{4.14} 
\begin{split}
& \cos^2(\alpha_0^{2,k_2})=|\zeta(w_2)|,\ w_2\in W_2, \\ 
& \cos(2\alpha_0^{3,k_3})=|\zeta(w_3)|,\ w_3\in W_3, 
\end{split}
\ee 
where 
\bdis 
\begin{split}
& w_2=w_2(U,\pi L,k_2,f_2,\sigma_0^2-\delta,\sigma_0^2+\delta), \\ 
& w_3=w_3(U,\pi L,k_3,f_3,\sigma_0^3-\delta,\sigma_0^3+\delta) 
\end{split}
\edis  
for every set of fixed and admissible corresponding parameters. 

\begin{remark}
Of course, 
\bdis 
\begin{split}
& S^l_{\sigma_0^l,\delta}\bigcap S^k_{\sigma_0^k,\delta}=\emptyset,\ k\not=l;\ k,l=1,2,3 
\Rightarrow \\ 
& W_l\bigcap W_k=\emptyset,\ l\not=k. 
\end{split} 
\edis 
\end{remark} 

\subsection{} 

Next, we choose an arbitrary finite set 
\bdis 
\{U_n\}:\ 0<U_1<U_2<\dots<U_{n_0}<\frac{\pi}{4},\ n_0\in\mbb{N}. 
\edis  
Every fixed $U_n$ generates the following three infinite sets 
\be \label{4.15} 
W_1(U_n),\ W_2(U_n),\ W_3(U_n),\ n=1,\dots,n_0. 
\ee  
Now, for every fixed $n$ we choose a single element 
\be \label{4.16} 
w_1^n,w_2^n,w_3^n,\ n=1,\dots,n_0
\ee  
from every corresponding set in (\ref{4.15}). Consequently, we have (comp. (\ref{4.13}), (\ref{4.14})) the following grafts (=equalities) 
\be \label{4.17} 
\begin{split}
& \sin^2(\alpha_0^{1,k_1,n})=|\zeta(w_1^n)|, \\ 
& \cos^2(\alpha_0^{2,k_2,n})=|\zeta(w_2^n)|, \\ 
& \cos(2\alpha_0^{3,k_3,n})=|\zeta(w_3^n)|, \\ 
& n=1,\dots,n_0 
\end{split}
\ee  
where (comp. (\ref{4.4})) 
\bdis 
\alpha_0^{1,k_1,n}=\alpha_0(U_n,\pi L,k_1;f_1),\dots 
\edis 

\subsection{} 

Finally, we make a grafting (we use the equalities (\ref{4.17}) as substitutions) on the exact complete hybrid formula (\ref{3.6}) (=our tree). As the result of these operations we obtain the following 

\begin{mydef11}
For every finite set 
\bdis 
\{U_n\}:\ 0<U_1<U_2<\dots<U_{n_0}<\frac{\pi}{4}
\edis  
there are the following elements 
\bdis 
w_1^n\in S^1_{\sigma_0^1,\delta},\ 
w_2^n\in S^2_{\sigma_0^2,\delta},\ 
w_3^n\in S^3_{\sigma_0^3,\delta},\ 
n=1,\dots,n_0
\edis 
where 
\be \label{4.18} 
\begin{split}
& w_1^n=w_1(U_n,\pi L,k_1,f_1,\sigma_0^1-\delta,\sigma_0^1+\delta), \\ 
& w_2^n=w_2(U_n,\pi L,k_2,f_2,\sigma_0^2-\delta,\sigma_0^2+\delta), \\ 
& w_3^n=w_3(U_n,\pi L,k_3,f_3,\sigma_0^3-\delta,\sigma_0^3+\delta), 
\end{split}
\ee 
such that the following exact meta-functional equation holds true: 
\be \label{4.19} 
\begin{split}
& \left\{
\prod_{r=1}^{k_2}\frac{\tilde{Z}^2(\alpha_r^{2,k_2,n})}{\tilde{Z}^2(\beta_r^{k_2,n})}
\right\}|\zeta(w_2^n)|- 
\left\{
\prod_{r=1}^{k_1}\frac{\tilde{Z}^2(\alpha_r^{1,k_1,n})}{\tilde{Z}^2(\beta_r^{k_1,n})}
\right\}|\zeta(w_1^n)|= \\ 
& =
\left\{
\prod_{r=1}^{k_3}\frac{\tilde{Z}^2(\alpha_r^{3,k_3,n})}{\tilde{Z}^2(\beta_r^{k_3,n})}
\right\}|\zeta(w_3^n)|, \\ 
& \forall\- L\geq L_0,\ 1\leq k_1,k_2,k_3\leq k_0, 
\end{split}
\ee  
where 
\be\label{4.20} 
\begin{split}
& \alpha_r^{1,k_1,n}=\alpha_r(U_n,\pi L,k_1;f_1), \\ 
& \beta_r^{k_1,n}=\beta_r(U_n,\pi L,k_1),\\ 
& \dots , \\ 
& n=1,\dots,n_0. 
\end{split}
\ee 
\end{mydef11}

\subsection{} 

Next, we obtain by the method (\ref{4.1}) -- (\ref{4.17}) other new-type formula. 

\begin{mydef12}
Under the same assumptions as in our Theorem 1 we have the following exact meta-functional equation: (\ref{3.2}) $\Rightarrow$ 
\be \label{4.21} 
\begin{split}
	& \left\{
	\prod_{r=1}^{k_1}\frac{\tilde{Z}^2(\alpha_r^{1,k_1,n})}{\tilde{Z}^2(\beta_r^{k_1,n})}
	\right\}|\zeta(w_1^n)|+ 
	\left\{
	\prod_{r=1}^{k_2}\frac{\tilde{Z}^2(\alpha_r^{2,k_2,n})}{\tilde{Z}^2(\beta_r^{k_2,n})}
	\right\}|\zeta(w_1^n)|=1 \\  
	& \forall\- L\geq L_0,\ 1\leq k_1,k_2\leq k_0. 
\end{split}
\ee  
\end{mydef12} 

\section{Asymptotic meta-functional equations} 

\subsection{} 

First, the method (\ref{4.1}) -- (\ref{4.17}) applied on the asymptotic formula (\ref{3.3}) gives the following 

\begin{mydef41}
The following asymptotic meta-functional equation 
\be \label{5.1} 
\begin{split}
& \left\{
\prod_{r=1}^{k_1}\frac{\left|\zeta(\frac 12+i\alpha_r^{1,k_1,n})\right|^2}{\left|\zeta(\frac 12+i\beta_r^{k_1,n})\right|^2}
\right\}|\zeta(w_1^n)|+
\left\{
\prod_{r=1}^{k_2}\frac{\left|\zeta(\frac 12+i\alpha_r^{2,k_2,n})\right|^2}{\left|\zeta(\frac 12+i\beta_r^{k_2,n})\right|^2}
\right\}|\zeta(w_2^n)|\sim \\ 
& \sim 1,\ L\to\infty  
\end{split} 
\ee 
holds true. 
\end{mydef41} 

\subsection{} 

Second, the method (\ref{4.1}) -- (\ref{4.17}) applied on the asymptotic formula (\ref{3.7}) gives the following 

\begin{mydef42}
	\be \label{5.2} 
	\begin{split}
		& \left\{
		\prod_{r=1}^{k_2}\frac{\left|\zeta(\frac 12+i\alpha_r^{2,k_2,n})\right|^2}{\left|\zeta(\frac 12+i\beta_r^{k_2,n})\right|^2}
		\right\}|\zeta(w_2^n)|-
		\left\{
		\prod_{r=1}^{k_1}\frac{\left|\zeta(\frac 12+i\alpha_r^{1,k_1,n})\right|^2}{\left|\zeta(\frac 12+i\beta_r^{k_1,n})\right|^2}
		\right\}|\zeta(w_1^n)|\sim \\ 
		& \sim 
		\left\{\prod_{r=1}^{k_3}\frac{\left|\zeta(\frac 12+i\alpha_r^{3,k_3,n})\right|^2}{\left|\zeta(\frac 12+i\beta_r^{k_3,n})\right|^2}
		\right\}|\zeta(w_3^n)|,\\ 
		& L\to\infty  . 
	\end{split} 
	\ee 
\end{mydef42} 

\subsection{} 

Now we obtain (see (\ref{5.1}), (\ref{5.2})) the following 

\begin{mydef43}
\be \label{5.3} 
\begin{split}
&  \left\{
\prod_{r=1}^{k_2}\frac{\left|\zeta(\frac 12+i\alpha_r^{2,k_2,n})\right|^2}{\left|\zeta(\frac 12+i\beta_r^{k_2,n})\right|^2}
\right\}|\zeta(w_2^n)|\sim \\ 
& \sim \frac 12+\frac 12
\left\{\prod_{r=1}^{k_3}\frac{\left|\zeta(\frac 12+i\alpha_r^{3,k_3,n})\right|^2}{\left|\zeta(\frac 12+i\beta_r^{k_3,n})\right|^2}
\right\}|\zeta(w_3^n)|,\ L\to\infty, 
\end{split}
\ee 
\be \label{5.4} 
\begin{split}
& \left\{
\prod_{r=1}^{k_1}\frac{\left|\zeta(\frac 12+i\alpha_r^{1,k_1,n})\right|^2}{\left|\zeta(\frac 12+i\beta_r^{k_1,n})\right|^2}
\right\}|\zeta(w_1^n)|\sim \\ 
& \sim \frac 12-\frac 12
\left\{\prod_{r=1}^{k_3}\frac{\left|\zeta(\frac 12+i\alpha_r^{3,k_3,n})\right|^2}{\left|\zeta(\frac 12+i\beta_r^{k_3,n})\right|^2}
\right\}|\zeta(w_3^n)|,\ L\to\infty. 
\end{split}
\ee 
\end{mydef43} 

\begin{remark}
The set of asymptotic meta-functional equations (\ref{5.1}) -- (\ref{5.4}) represents the higher level of $\zeta$-encoding (in comparison with this one that represents the corresponding set of asymptotic complete hybrid formulas) of the set of elementary trigonometric formulas 
\bdis 
\begin{split}
& \sin^2x+\cos^2x=1,\ \cos^2x-\sin^2x=\cos 2x, \\ 
& \cos^2x=\frac 12+\frac 12\cos 2x,\ \sin^2x=\frac 12-\frac 12\cos 2x. 
\end{split}
\edis 
\end{remark}

\subsection{} 

Putting in (\ref{5.2}) 
\be \label{5.5} 
k_1=k_2=k_2=k;\ 1\leq k\leq k_0, 
\ee
we obtain the following 

\begin{mydef44}
\be \label{5.6} 
\begin{split}
& \left\{\prod_{r=1}^{k}\left|\zeta(\frac 12+i\alpha_r^{2,k,n})\right|^2\right\}|\zeta(w_2^n)|-
\left\{\prod_{r=1}^{k}\left|\zeta(\frac 12+i\alpha_r^{1,k,n})\right|^2\right\}|\zeta(w_1^n)|\sim \\ 
& \sim 
\left\{\prod_{r=1}^{k}\left|\zeta(\frac 12+i\alpha_r^{3,k,n})\right|^2\right\}|\zeta(w_3^n)|,\ L\to\infty . 
\end{split}
\ee 
\end{mydef44} 

\begin{remark}
We shall call the formula (\ref{5.6}) as the secondary asymptotic meta-functional equation since it is the result of immediate metamorphosis of the formula (\ref{5.2}) in the case (\ref{5.5}). 
\end{remark} 

\begin{remark}
The case $k=1$ in (\ref{5.5}) gives the formula whose simple variant is the formula (\ref{1.2}) in the Introduction. 
\end{remark} 

\section{Synergetic formulas only for Riemann's zeta-function: interactions between some parts of this function} 

\subsection{} 

Let us remind that every element from the complete hybrid formulas set (\ref{3.2}), (\ref{3.3}), (\ref{3.5}) -- (\ref{3.7}) is synergetic one - see our interpretation of this phenomenon in our paper \cite{8}. Really, every from these formulas is the result of an interaction controlled by the Jacob's ladder $\vp_1(t)$ between following sets 
\bdis 
\left\{\left|\zf\right|^2\right\},\ \{\sin^2t\},\ \{\cos^2t\},\ \{\cos 2t\},\ t\geq L_0. 
\edis 

\subsection{} 

In this paper, we have introduced some interactions of Riemann's zeta-function with itself on the critical strip. Namely, the interactions between following four parts 
\be \label{6.1} 
\left\{\left|\zf\right|^2\right\},\ t\geq L_0,\ |\zeta(s)|,\ s\in S^l_{\sigma_0^l,\delta},\ l=1,2,3
\ee  
of the Riemann's zeta-function 
\bdis 
\zeta(s),\ s\in\mbb{C}\setminus \{1\}. 
\edis 
Here, of course, 
\bdis 
S^l_{\sigma_0^l,\delta}\subset (\sigma_1,\sigma_2)\times (0,+\infty),\ 
S^l_{\sigma_0^l,\delta}\bigcap S^k_{\sigma_0^k,\delta}=\emptyset,\ k\not=l,\ k,l=1,2,3. 
\edis 
Let us call the set of these interactions as the $\zeta$-chemical reaction of the \emph{substances} (\ref{6.1}) - this is some analogue to what we have based on the classical Belousov-Zhabotinski chemical reactions. 

\begin{remark}
The result of these $\zeta$-chemical reactions (=$\zeta$-chemical compounds) 
are our synergetic meta-functional equations (\ref{4.19}), (\ref{4.21}) as well as the formulae (\ref{5.1}) -- (\ref{5.4}), (\ref{5.6}). 
\end{remark} 

\begin{remark}
Since the Planck length reads 
\bdis 
L_P=8.1\times 10^{-35}cm
\edis  
and the Planck time reads 
\bdis 
T_P=2.7\times 10^{-43}s 
\edis 
then, from a point of view of the sciences, we may use (for example) the following complementary condition on the set 
\bdis 
\{ U_n\}:\ U_{n+1}-U_n>10^{-43},\ n=1,\dots,n_0-1, U_1,\frac{\pi}{4}-U_{n_0}>10^{-43}. 
\edis  
This condition corresponds with the point of view of Jakov Zeldovich for using maths on study of real-world phenomena. 
\end{remark}

I would like to thank Michal Demetrian for his moral support of my study of Jacob's ladders.

\end{document}